\documentclass{amsart}
\usepackage{amsfonts,amssymb,amscd,amsmath,enumerate,verbatim,calc}
\usepackage[all]{xy}
\SelectTips{eu}{}

\newcommand{\bsx}{\boldsymbol{x}}
\newcommand{\bsy}{\boldsymbol{y}}

\newcommand{\mf}{\mathfrak }

\newcommand{\mit}{\mathit}

\newcommand{\BN}{\mathbb N}
\newcommand{\BZ}{\mathbb Z}

\newcommand{\lra}{\longrightarrow}
\newcommand{\xra}{\xrightarrow}
\newcommand{\arr}{\rightarrow}

\newcommand{\Image}{\operatorname{Im}}
\newcommand{\Homology}{\operatorname{H}}

\newcommand{\grade}{\operatorname{grade}}
\newcommand{\depth}{\operatorname{depth}}

\newcommand{\cmd}{\operatorname{cmd}}

\newcommand{\Hom}{\operatorname{Hom}}
\newcommand{\hh}{\operatorname{H}}

\newcommand{\length}{{\ell}}

\theoremstyle{plain}

\theoremstyle{definition}

\theoremstyle{plain}

\newtheorem{theorem}{Theorem}[section]

\newtheorem{lemma}[theorem]{Lemma}

\newtheorem*{maintheorem}{Main Theorem}

\theoremstyle{definition}

\newtheorem{example}[theorem]{Example}

\newenvironment{thmprooflc}{%

\begin{proof}}{\end{proof}}

\newenvironment{thmproofpoin}{%

\begin{proof}}{\end{proof}}

\newtheorem{chunk}[theorem]{}

\theoremstyle{remark}

\newtheorem{remark}[theorem]{Remark}

\numberwithin{equation}{theorem}

\input xy
\xyoption{all}


\begin{document}

\title[ Free resolutions of parameter ideals]
{Free Resolutions of parameter ideals for some rings with finite
local cohomology}

\author{Hamid Rahmati}
\address{Department of Mathematics, University of Nebraska,  Lincoln, NE 68588, U.S.A.}
\email{hrahmati@math.unl.edu}

\begin{abstract}
Let $R$ be a $d$-dimensional  local ring, with maximal ideal $\mf m$, containing a field  and let $ x_1,
\dots , x_d$ be a system of parameters for $R$. If $\depth R \geq
d - 1$ and the local cohomology module $\hh_{\mit
m}^{d-1}(R)$ is finitely generated, then there exists an integer $n$
such that  the modules $R/(x_1^i, \dots, x_d^i)$
have the same Betti numbers, for all $i \geq n$.

\end{abstract}

 \maketitle

\section*{Introduction}

 Throughout this note $R$ is a local noetherian ring with
maximal ideal $\mf m$.

Let $M$ be a finitely generated $R$-module and let $F$ be a minimal free resolution of $M$.
The Poincar\'e series of $M$ is the formal power series
$\operatorname P_M^R(t) =  \sum_{i=0}^\infty (\operatorname {rank}
F_i)t^i$. We  let
$\Omega^i_R(M)$ denote the $i$th syzygy of $M$, that is to say $\operatorname
{Coker}
\partial^F_{i+1}$.
Let  $\bsx=x_1, \dots ,x_d$ be a system
 of parameters for $R$. For each $n \in  \BN$, let $\bsx^n$ denote the sequence
$x_1^n,\dots,x_d^n$.

Several  classical results in local algebra establish that ideals contained in
the large powers of the maximal ideal exhibit  a similar behavior.
In his thesis, \cite {Lai}, Y. H. Lai considers the following question, which he attributes  to D. Katz: Do the
Poincar\'e series of the modules $R/(\bsx^n)$ behave uniformly, for
$n$ large enough?

Clearly, the answer is  yes if $R$ is Cohen-Macaulay: In this
case $\bsx^n$ is a regular sequence, for each $n \geq 1$, so
$R/(\bsx^n)$ is resolved by a Koszul complex  on $d$ elements. It is also
not difficult to obtain a positive
 answer when $\dim R =1$.

 Lai proves that if $\dim R = 2$, $\depth R = 1$ and the
local cohomology module $\hh_\mf m^1(R)$ is finitely generated, then for $n$ large enough one has
\[
\operatorname P^R_{R/(\bsx^n)}(t) = 1+2t+t^2 + t^2\operatorname P^R_{H^1_{\mf m}(R)}(t).
\]

Now we state our main result. It evidently covers the first two cases mentioned above.
Also part (ii) of our main theorem generalizes Lai's result, since the Canonical Element Conjecture holds for  2-dimensional rings.

\begin{maintheorem}
\label {m1} Let $R$ be $d$-dimensional local ring with maximal ideal $\mf m$. If $d- \depth R \leq 1$ and the $R$-module $H =
\Homology _{\mf m}^{d-1}(R)$ is  finitely generated, then there
exists an integer $n$, such that for each system of parameters $\bsx$
for $R$ contained in $\mf m^n$ the following assertions hold

\begin{enumerate}[\quad\rm(i)]
\item  $\Omega^{d+1}_R(R/(\bsx)) \cong
\Omega^{d-1}_R(H)$.

\item If in addition, the Canonical Element Conjecture holds for $R$,
then one has
\[
\operatorname P^R_{R/(\bsx)}(t) = (1+t)^d + t^2\operatorname P^R_{H}(t).
\]
\end{enumerate}
\end{maintheorem}

The \emph{ Cohen-Macaulay defect} of $R$ is the number $\cmd R =\dim R - \depth R $. One always has
$\cmd R \geq 0$, and equality  characterizes Cohen-Macaulay rings. On the other hand,
Cohen-Macaulay rings are also characterized by the equality $\Homology _{\mf m}^i(R)=0$ for all $i\ne d$.
Comparing these conditions with
the hypotheses on $R$ in the theorem, one may say that these hypotheses, in some sense, give  the least possible
extension of Cohen-Macaulay rings.

\section{Free resolutions of almost complete intersection ideals}
Let $\boldsymbol y = y_1,\dots,y_n$ be a sequence in $R$. We let
$K(\boldsymbol y;M)$ denote the Koszul complex on $\boldsymbol y$
with coefficients in $M$. Set
\[
 \Homology_i(\boldsymbol y;M) =
\hh_i(K(\boldsymbol y;M)).
\]

If  $P(t)=\sum_{i=0}^\infty a_i t^i$  and    $Q(t)=\sum_{i=0}^\infty b_i t^i$  are formal power, we write
$P \preccurlyeq Q$ to indicate that $a_i\leq b_i$ holds for all $i\geq 0$.


The main result of this section is the following theorem:

\begin{theorem}
\label{poincare} Let $R$ be a $d$-dimensional ring and let $\bsx$ be a system of parameters for $R$.
Set $H_i=\hh_i(\bsx;R)$ for $i=1,\dots,d$. One then has
\[
\operatorname P^R_{R/(\bsx)}(t) \preccurlyeq (1+t)^d+ \sum_{i=1}^{\cmd R }
t^{i+1}\operatorname P^R_{H_i}(t).
\]
Moreover, if  $\cmd R \leq 1 $ and the Canonical Element Conjecture holds for
$R$, then one has
\[
\operatorname P^R_{R/(\bsx)}(t) =(1+t)^d+
t^2\operatorname P^R_{H_1}(t).
\]
\end{theorem}
Under the hypotheses of the second part of the theorem, $\bsx$ generates an almost complete intersection ideal.
 Some of the discussion below is carried out in this more general framework.
 Recall that an ideal $I$ is called  \emph{almost complete intersection}
if $\operatorname {grade} I \geq \mu(I) -1$, where $\operatorname {grade} I$ is the  maximal length of an
$R$-regular sequence in $I$ and $\mu(I)$ is the number of minimal generators of $I$.

Now we give an example that shows the inequality in Theorem~\ref{poincare} can be strict if $\depth R < d-1$.
\begin{example} \label{example}
Set $R=k[\![a,b,c]\!]/(ac,bc,c^2)$, then $\dim R=2$ and $\depth R =0$. Consider the system of parameters $\bsx = a,b$. Using Macaulay 2 we get:
\begin{align*}
P^R_{H_2}(t) & = 1+3t+6t^2+13t^3+28t^4+\cdots \\
P^R_{H_1}(t) & = 3+7t+12t^2+26t^3+56t^4+\cdots \\
P^R_{R/(\bsx)}(t) & = 1+2t+3t^2+7t^3+15t^4+\cdots.
\end{align*}
Thus one has
\begin{align*}
P^R_{R/(\bsx)}(t)\prec  1+2t+4t^2+8t^3+15t^4+\cdots \\
= (1+t)^2 + t^2 P^R_{H_1}(t) + t^3P^R_{H_2}(t) .
\end{align*}

\end{example}

Let $X$ be a complex of $R$-modules and let  $\partial^X$ denote its differential, set
\[
\operatorname {sup} \hh(X) = \operatorname {sup} \{n \in \BZ \mid
\hh_n(X) \ne 0 \}.
\]
Let $\Sigma^t$ denote the shift functor defined by
\[
(\Sigma^t X)_n = X_{n-t} \quad \mbox{ and } \quad  \partial^{\Sigma
X}_n = (-1)^t \partial^X_{n-t}.
\]
Let  $\alpha \colon X \arr Y$ be a morphism of complexes. Recall that
mapping cone of $\alpha$ is defined to be the complex $C$ such that
$C_n = X_{n-1} \oplus Y_n$ and \[   \partial^C_n ((x,y)) =
(-\partial^X_{n-1}(x),\partial^Y_n(y) + \alpha_{n-1}(x))
 \quad \mbox{for all} \quad (x,y) \in C_n.
\]
 A quasi-isomorphism is a morphism of complexes that induces isomorphism in homology in all degrees.

\begin{lemma}
\label{cone}
 Let $X$ be a complex with $s=\operatorname
{sup}\Homology (X) < \infty$ and let $F$ be a free resolution of
$\Homology_s(X)$. There exists a morphism of complexes $\alpha \colon
\Sigma^s F \arr X$ such that the mapping cone $X'$ of $\alpha$ satisfies
\begin{equation*}
\Homology_i(X') \cong
\begin{cases}
0 & \quad  \text {$i \geq s$} \\
 \Homology_i(X) & \quad  \text {$ i \leq s-1$}.
\end{cases}
\end{equation*}
\end{lemma}

\begin{proof}
Let $\tau_{\geq s}(X)$ be the complex
\[
 \cdots \lra X_{s+2} \lra X_{s+1} \lra \operatorname {ker} \partial^X_s
\lra 0
\]
and $\iota \colon \tau_{\geq s}(X) \to X$ be the inclusion map.
Since $F$ is a bounded below complex of free modules,
 there is a morphism of complexes $\beta \colon F \arr \tau_{\geq s}(X)$ such that the
following diagram is commutative

\begin{equation*}
\xymatrix{ \Sigma^sF \ar[r]^\beta \ar[dr]_{\Sigma^s\epsilon} &
\tau_{\geq s}(X) \ar[r]^{\hskip 12pt \iota} \ar[d]^\pi & X \\ &
\Sigma^s\Homology_s(X), }
\end{equation*}
where $\epsilon \colon F \to \hh_s(X) $  and $\pi$ are
quasiisomorphism. Set $\alpha = \iota
\beta$ and let $C$ be the mapping cone of $\alpha$. One has an exact sequence
\[
0 \lra X \lra X' \lra \Sigma^{s+1} F \lra 0
\]
of complexes.
It induces an exact sequence of homology modules
\begin{align*}
\cdots \arr \Homology_{i+1}(X) \arr \Homology_{i+1}(X') \arr
 \Homology_{i+1}(\Sigma^{s+1}F) \xra {\Homology_i(\alpha)} \Homology_i(X) \arr
 \cdots.
\end{align*}
From the construction of $\alpha$ one sees that $\hh_s(\alpha)$
is an isomorphism, and since $\hh_i(\Sigma^{s+1}F) \cong
\hh_{i-s-1}(F) = 0$ for all $i \ne s+1$, we get the desired
result.
\end{proof}

\begin{remark}\label{nonminimality}
Assume $X$ is a complex of free modules such that $X_i=0$, for $i<0$, and $\Homology_i(X)=0$, for $i >s$,
where $s$ is some positive integer. Let $F^i$  be a free resolution of $\Homology_i(X)$, for $i=1,\dots,s$.
Applying Lemma~\ref{cone} $s$ times, one gets
a free resolution $G$ of $\Homology_0(X)$ such that $G_i=X_i \oplus  F^1_{i-2} \oplus \cdots \oplus F^s_{i-s-1}$.
However, even if the complexes $X,F^1,\dots,F^s$ are minimal, $G$ need not be minimal; see Example~\ref{example}.
\end{remark}

Now we show the relation between the Poincar\'e series of an almost
complete intersection ideal and the Poincar\'e series of  its first Koszul homology module.
\begin{lemma}
\label {syzygy} Let $I$ be an almost intersection ideal of $R$ and
let $\boldsymbol y = y_1, \dots, y_r $ be a  minimal set generators for $I$.  Set
$H=\hh_1(\boldsymbol y;R)$. One then has
\begin{align*}
 \Omega^{r+1}_R(R/I) & =  \Omega^{r-1}_R(H)  \\
 \operatorname P^R_{R/I}(t) & =  Q(t) + t^{r+1}\operatorname
P^R_{N}(t),
\end{align*}
 where $Q(t)$ is a polynomial of degree $r$ and $N =
\Omega^{r-1}_R(H)$.

\end{lemma}

\begin{proof}
Let $K$ denote the Koszul complex of $\bsy$ with coefficients in
$R$. Since $\grade I \geq r-1$, one has $\hh_i(\bsy;R) = 0$ for $i \ne
0,1$. Let $F $ be a minimal free resolution of $H$. Let $\alpha
\colon \Sigma F \arr K$ be the map from Lemma~\ref {cone} and let
$C$ be the mapping cone of $\alpha$. The complex $C$ is a complex of
free modules and has only one nonvanishing homology module namely,
$\hh_0(C)\cong \hh_0(\bsy;R) \cong R/I$, so $C$ is a free
resolution of $R/I$.

By construction of mapping cones the complex $C$ is minimal if and
only if $\alpha F \subseteq \mf mK$. Since $F$ is minimal and $C_m=F_m$ for all $m \geq r+1$,
non-minimality can only happen in the first $r+1$ degrees, therefore
\[
\Omega^{r+1}_R(R/I) = \Omega^r_R(H).
\]
This completes the proof of the first equality. The second
equality follows from the first one.
\end{proof}

\begin{chunk}
For a system of parameters $\bsx$ for $R$, let $F_{\bsx} $ be a free
resolution of $R/(\bsx)$ and let $\gamma_{\bsx} \colon K_{\bsx} =
K(\bsx;R) \to F_{\bsx}$ be a lifting of the map
$ K_{\bsx} \to R/(\bsx)$. The
following are equivalent:

\begin{enumerate}[\quad\rm(i)]
\label{cec}
\item The Canonical Element Conjecture holds for $R$.
\item For every system of parameters $\bsx$ for $R$ and every free resolution
$F_{\bsx}$ of $R/(\bsx)$, the map $\hh(k\otimes_R \gamma_{\bsx}) \colon \hh(k\otimes_R K_{\bsx}) \to \hh(k \otimes_R F_{\bsx})$ is injective.

\end{enumerate}

It is shown in \cite [(1.6)]{AI} that the equivalence of (i) and
(ii) follows from a theorem of P. Roberts \cite{Ro}. A proof of
Roberts' theorem is given in \cite [(1.3)]{HK}.

The Canonical Element Conjecture holds for $R$
provided $R$ is equicharacteristic or $\dim R \leq 3$:
 Hochster has proved that the Canonical Element Conjecture is equivalent to the Direct Summand Conjecture,
 and  that the conjectures hold if $R$ is equicharacteristic or $\dim R \leq 2$, see \cite{Hoch}. In \cite{He},
 R. Heitmann shows that the Direct Summand Conjecture, hence the Canonical Element Conjecture, holds for every  3-dimensional ring.

\end{chunk}

\begin{thmproofpoin}
 The inequality follows from Remark~\ref{nonminimality}.

  For the rest of the proof, we keep the notation in the proof of Lemma \ref{syzygy}. To prove the equality it
suffices to show that $C$ is minimal. Assume not, and let $0 \leq n \leq
d-1$ be such that $\partial (C_n) \nsubseteq \mf mC_{n-1}$. Since
$\Image (\partial^K)$ and $\Image (\partial^F)$ are in $\mf mC$, there
exists an element $f$ such that $\partial (f)=e \in K_n\backslash \mf mK_n$.
It follows that $f$ is not in $  \mf mF_n$ hence $Rf$ is a direct summand of $C_n$,
 and $Re$ is a direct summand of $K_n$. Let $G$ be
the complex $0 \arr Rf  \xra {\lambda} Re \arr 0$ where $Re$ is in
degree $n$ and $\lambda$ is the restriction  of $\partial^{C}_n$.  Set
$\overline{C} = C/G$. Since $\lambda$ is an isomorphism, $\overline
X'$ is exact, hence is a free resolution of $R/(\bsx)$. Let $\pi \colon
C \arr \overline C$ be the natural surjection.

The inclusion map
$\iota \colon K \arr C$ is a lifting  of the augmentation map
$\alpha \colon K \arr R/(\bsx)$. Thus $\pi \iota$ is a lifting of
$\alpha$, but $\hh_n(k \otimes_R \pi \iota)(1 \otimes e) = 0$. This
contradicts, by \ref {cec}, the hypothesis that Canonical Element
Conjecture holds for $R$. Therefore, $C$ is minimal; this implies
the theorem.
\end{thmproofpoin}

\section{standard system of parameters}

Let $M$ be a finitely generated $R$-module. A system of parameters
$\bsx$ for $M$ is said to be \emph {standard} if
\[
(\bsx)\Homology ^i_{\mf m}(M/{(x_1,\dots,x_j)M}) =0
\]
 holds, for all
non-negative integers $i,j$ with $i + j < d$. For information about
standard system of parameters we refer the reader to \cite {SV} and
\cite {Tr}.

An $R$-module $M$ is said to have \emph{finite local cohomology} if
for each integer $i \leq \dim M-1$ the local cohomology module
$\Homology^i_{\mf m}(M)$ is of finite length. Modules with finite
local cohomology are also called \emph {generalized Cohen-Macaulay  modules}.

The following statement is a consequence of \cite [ (2.1) and
(3.1)]{Tr} and \cite [(3.7)]{STC}.

\begin{chunk}
\label{Trung1}
 An $R$-module $M$ has finite local cohomology if and only if
 there exists a positive integer $n$ such that every system of parameters in
 $\mf m^n$  is standard.
\end{chunk}
For an $R$-module $M$ we set
$ \hh^i(\boldsymbol
y;M) = \hh_{-i}(\Hom_R(K(\boldsymbol y ; R),M))$.
One then has $\hh_i(\boldsymbol y;M) \cong \hh^{d-i}(\boldsymbol
y;M)$, for all $i$, see \cite [1.6.10]{BH}.

\begin{theorem}
\label{lc}
 Let $M$ have finite local cohomology. Set $g = \depth M$. If
$\bsx$ is a standard system of parameters for $M$,  then for every
integer $n \geq 1$, the canonical map $\lambda_n
\colon\Homology^g(\bsx^n;M) \arr \Homology^g_{\mf m}(M)$ is an
isomorphism.
\end{theorem}

To prove this theorem we need to recall some facts about standard system
of parameters and modules with finite local cohomology.

The next statement follows from \cite [(3.1)]{Tr} and
\cite [(3.3)]{STC}.

\begin{chunk}
\label{regularity}Let $\bsx$ be a system of parameters for $M$. If
$M$ has finite local cohomology and $\depth M =g$, then the
subsequence $x_1,\dots , x_g$ of $\bsx$ is $M$-regular.
\end{chunk}

In  \cite [(1)]{Ho}, standard systems of parameters are
characterized in terms of Koszul homology:
\begin{chunk}
 \label{Hoa1}
Assume $M$ has finite local cohomology. Let $\bsx$ be a system of
parameters for $M$. The following are then equivalent:
\begin{enumerate}[\quad\rm(i)]

\item $\bsx$ is standard.

\item $\length (\Homology_p(x_1, \dots,x_r;M)) = \length
(\Homology_p(x_1^2,\dots,x_r^2;M))$, for $p\geq 1$  and $1\leq r \leq d$.

\item $\length ( \Homology_p(x_1,\dots,x_r ; M)) =
\sum_{i=0}^{r-p} \binom r{i+p} \length(\Homology^i_{\mf m}(M))$, for
all  $p\geq 1$ and $1\leq r \leq d$.
\end{enumerate}
In particular, if $\bsx$ is standard, then \[\length (
\Homology_{d-g}(\bsx ; M)) = \length(\Homology^g_{\mf m}(M)),\quad
\mbox {where} \quad g=\depth M \quad \mbox{and} \quad d=\dim M.\]
Here $\length(M)$ denotes the length
of $M$.
\end{chunk}

The next result is \cite [(4)]{Ho}.

\begin{chunk}

\label{inequality}
 Let $M$  have finite local cohomology, and let $x_1,\dots, x_d$ be a
system of parameters for $M$. Let $n_1, \dots, n_d$ be positive integers. For all
$p>0$ and  $1\leq r \leq d$ and for all positive integers $m_1,
\dots ,m_r$ satisfying $n_1 \leq m_1, \dots, n_r \leq m_r$, one has
\[
\length(\Homology_p(x_1^{n_1}, \dots , x_r^{n_r};M)) \leq
\length(\Homology_p(x_1^{m_1}, \dots , x_r^{m_1};M)).
\]
\end{chunk}

\begin{lemma}
\label{equality}
 Let $M$ have finite local cohomology
and let $\bsx$ be a standard system of parameters for M. For all positive integers $p$ and $m$, one has

\[
\length (\Homology_p(\bsx;M)) = \length (\Homology_p(\bsx^m;M)).
\]
\end{lemma}

\begin{proof}

For every $i > 0$, from \ref{Hoa1}, one obtains:
\[
\length (\Homology_p(\bsx ;M))=\length (\Homology_p(\bsx^{2i};M)).
\]
\par For each integer $m > 0$ there exists an integer $i>0$ such
that $ 2i \geq m$. Using \ref{inequality}, we get inequalities
\begin{eqnarray*}
 \length (\Homology_p(\bsx;M)) \leq   \length
(\Homology_p(\bsx^m;M)  \leq \length (\Homology_p(\bsx^{2i};M)),
\end{eqnarray*}
which imply the desired statement.
\end{proof}

Next we recall the relation between Koszul homology and  local
cohomology . Let $\bsx$ be a system of parameters for $M$ and set
$K^{(n)}=K(\bsx^n;R)$ and let $e_1^{(n)},\dots,e_d^{(n)}$ denote the
standard basis of $K^{(n)} \simeq R^d$. For all $n \geq 1$, there is
a commutative diagram
\begin{equation*}
\xymatrix{ K_1^{(n+1)} \ar[d] \ar[r]^-{\varphi^{(n)}_1}&  K_1^{(n)}
\ar[d]  \\  K_0^{(n+1)} \ar@{=}[r] & K_0^{(n)} }
\end{equation*}
where $\varphi_1^{(n)}(e^{(n+1)}_j) = x_je^{(n)}_j$ for $j = 1,
\dots, d$. This defines a morphism of complexes $\varphi ^{(n)} =
\wedge \varphi^{(n)}_1 \colon K^{(n+1)} \to  K^{(n)}$. Set $\psi
^{(n)} = \Hom_R(\varphi ^{(n)} , M)$, then $\psi ^{(n)}$ induces a
map $\alpha _{(n)}^i \colon \hh^i(\bsx^n;M) \arr
\hh^i(\bsx^{(n+1)};M)$. By \cite[3.5.6]{BH}, one has
\[\Homology^{i}_{\mf m}(M)= \varinjlim \Homology^i(\bsx^{n};R)\quad \mbox{for
all} \quad i \geq 0.\]

\begin{lemma}
\label{commdiag} Let  $\boldsymbol z = z_1,\cdots , z_t$ be a
sequence in $R$ and let $s$ be the length of maximal $M$-regular
sequences in $(\boldsymbol z)$. For each positive integer $n$ , if
$\boldsymbol y = y_1, \cdots, y_r$ is an $M$-regular sequence in
$(\boldsymbol z^{n+1})$, then one has the following commutative diagram

\begin{equation*}
\xymatrix{
\hh^s(\boldsymbol z^n;M) \ar[r]^-{\alpha^s_{(n)}}   & \hh^s(\boldsymbol z^{n+1};M) \\
\hh^{s-r}(\boldsymbol z^n;\overline{M})
\ar[r]^-{\bar{\alpha}^{s-r}_{(n)}} \ar[u]_{\simeq}  &
\hh^{s-r}(\boldsymbol z^{n+1};\overline{M}) \ar[u]_{\simeq} }
\end{equation*}
where $\overline{M}= M/(\boldsymbol y)M$.

\end{lemma}

\begin{proof}
 If $r=1$, then one has a short exact sequence
\[
0 \lra M \xra {y_1} M \lra M/y_1M \lra 0.
\]
It induces a commutative diagram
\begin{equation*}
\xymatrix{ 0 \ar[r] & \hh^{s-1}(\boldsymbol z^n;M/y_1M)
\ar[r]^-{\eth_n^{s-1}} \ar[d]^{\alpha^{s-1}_{(n)}} &
\hh^{s}(\boldsymbol z^n;M) \ar[d]^{\alpha^s_{(n)}}
 \ar[r]^{y_1} & \hh^s(\boldsymbol z^n;M) \ar[d]^{\alpha^s_{(n)}} \\
0 \ar[r] & \hh^{s-1}(\boldsymbol z^{n+1};M/y_1M) \ar[r]
^-{\eth_{n+1}^{s-1}} & \hh^s(\boldsymbol z^{n+1};M)\ar[r]^{y_1}&
\hh^s(\boldsymbol z^{n+1};M).}
\end{equation*}
Since $y_1$ annihilates both $\hh^{g}(\boldsymbol z^n;M)$ and
$\hh^{g}(\boldsymbol z^{n+1};M)$, the connecting maps $\eth_n^{g-1}$ and
$\eth_{n+1}^{g-1}$ are isomorphisms and this gives  the desired
commutative diagram for $r=1$. The general case follows by iteration.
\end{proof}

\begin{thmprooflc}

Let $n \geq 1$. The sequence   $\bsx'=x_1, \dots ,x_g$ is
$M$-regular, see \ref {regularity}. So ${\bsx'}^{n+1}$ is
an $M$-regular sequence in the ideals $(\bsx^{n+1})  \subseteq (\bsx^n) $. Set
$\overline{M}=M/(\bsx')M$. By \ref {commdiag}, we get the following
commutative diagram
\begin{equation*}
\xymatrix{
\hh^g(\boldsymbol x^n;M) \ar[r]^-{\alpha^g_{(n)}}   & \hh^g(\boldsymbol x^{n+1};M) \\
\hh^0(\boldsymbol x^n;\overline{M})
\ar[r]^-{\bar{\alpha}^0_{(n)}} \ar[u]_{\simeq}  &
\hh^0(\boldsymbol x^{n+1};\overline{M}) \ar[u]_{\simeq}. }
\end{equation*}

The map $\bar{\alpha}^0_{(n)}$ is injective, since it  is
induced by the identity map in the following commutative diagram
\begin{equation*}
\xymatrix{
\Hom_R(K^{(n)},\overline{M})= \ar[d]^{\psi^{(n)}} &  0 \ar[r] & \overline{M} \ar[d]^= \ar[r] & \overline{M}^d \ar[d] \ar[r] &  \cdots \\
\Hom_R(K^{(n+1)},\overline{M})=  & 0 \ar[r] & \overline{M} \ar[r] &
\overline{M}^d \ar[r] & \cdots .}
\end{equation*}
 Therefore $\alpha ^g_{(n)}$ is injective.

  The modules $\Homology^g(\bsx^n;M)$ and $\Homology^g
(\bsx^{n+1};M)$ have the same length, see \ref{equality}, so  $\alpha^g_{(n)}$ is an
isomorphism and this implies the desired statement.

\end{thmprooflc}

\section{free resolutions of parameter ideals}
In this section we assume that R is a $d$-dimensional ring with finite local cohomology, and $\cmd R \leq 1$.

Because of the property recalled in \ref{Trung1}, the following result contains the Main Theorem stated in the introduction.

\begin{theorem}
Let $\bsx$ be a standard system of parameters for $R$ and set $H=\Homology_{\mf m}^{d-1}(R)$. One then has

\begin{enumerate}[\quad\rm(i)]
\item  $\Omega^{d+1}_R(R/(\bsx)) \cong
\Omega^{d-1}_R(H)$.

\item If in addition, the Canonical Element Conjecture holds for $R$,
then one has
\[
\operatorname P^R_{R/(\bsx)}(t) = (1+t)^d + t^2\operatorname P^R_{H}(t).
\]
\end{enumerate}
\end{theorem}

\begin{proof}

One has $\grade (\bsx) = \depth R$ and $\mu(\bsx) =\dim R$, so $ \grade (\bsx) \geq \mu (\bsx)-1$,
then part (i) follows from Lemma~\ref{syzygy} and Theorem~\ref{lc}, and part (ii) follows from Theorem~\ref{poincare} and Theorem~\ref{lc} .

\end{proof}

\section*{Acknowledgments}
I should like to express my gratitude to my advisers L. L. Avramov and S. Iyengar for their
support and guidance throughout the development of this paper.


\begin{thebibliography}{10}



\bibitem{AI}
L. L.~Avramov and S.~Iyengar, \emph{Gaps in Hochschild cohomology
imply smoothness for commutative algebra}, Math. Res. Letters
\textbf{12} (2005), 789-804.

\bibitem{BH}
W.~Bruns and J.~Herzog, \emph{Cohen-Macaulay rings}, vol.~39, Cambridge Stud.
  Adv. Math., Cambridge Univ. Press, 1998.

\bibitem{He}
R.~C. Heitmann, \emph{The direct summand conjecture in dimension 3}, Ann. of Math.
(2) \textbf{156} (2002), no. 2, 695-712.

\bibitem{Ho}
L.~T. Hoa, \emph{Koszul homology and generalized Cohen-Macaulay
modules}, Acta Math. Vietnam \textbf{18} (1993), no. 1, 91-98.

\bibitem{Hoch}
M.~Hochster, \emph{Canonical elements in local cohomology and the direct
summand conjecture}, J. Algebra \textbf{84} (1983), 503-553
\bibitem{HK}
C.~Huneke and J.~Koh, \emph{Some dimension $3$ cases of the
canonical element conjecture.},  Proc. Amer. Math. Soc. \textbf{98}
(1986), no. 3, 394-398.


\bibitem{Lai}
Y.~H. Lai, \emph{On the relation type of system of parameters and on
the Poincar\'e series of system of parameters.}, Ph.D Thesis Purdue
Unviversity, 1995

\bibitem{Ro}
P. ~Roberts, \emph{The equivalence of two forms of the Canonical
Element Conjecture}, undated manuscript.

\bibitem{STC}
P.~Schenzel, N.~V. Trung and N.~T. Coung, \emph{Verallgemeinerte  Cohen-Macaulay Moduln},
Math Nachr, \textbf{85} (1978), 57-73.

\bibitem{SV}
J.~St\"uckard and W.~Vogel, \emph{Buchsbaum rings and applications},
Berlin-Heidelberg-New York, Springer-Verlog, 1986.



\bibitem{Tr}
N.~V. Trung, \emph{Toward a theory of generalized Cohen-Macaulay
modules}, Nagoya Math. J. \textbf{102} (1986), 1-49.


\end{thebibliography}
\end{document}